\documentclass[11pt]{amsart}

\usepackage[utf8x]{inputenc} 
\usepackage[a4paper]{geometry} 
\usepackage{enumerate} 
\usepackage{amsmath} 
\usepackage{amssymb}
\usepackage{color,colortbl}
\usepackage[dvipsnames,table]{xcolor}
\usepackage[most]{tcolorbox}
\usepackage{fancyvrb}   
 
\definecolor{airforceblue}{rgb}{0.36, 0.54, 0.66}
\definecolor{bleudefrance}{rgb}{0.19, 0.55, 0.91}
\definecolor{darkorchid}{rgb}{0.6, 0.2, 0.8}
\definecolor{darkorange}{rgb}{1.0, 0.55, 0.0}
\definecolor{darkspringgreen}{rgb}{0.09, 0.45, 0.27}
\definecolor{commentoutput}{rgb}{0.40, 0.00, 0.0}
\definecolor{output}{rgb}{0.8, 0.0, 0.0}  
\definecolor{circOut}{rgb}{0.4, 1.0, 0.0} 

\definecolor{Gray}{gray}{0.7}  

\usepackage[arrow,curve,matrix,tips,frame,color]{xy}
\usepackage{adjustbox}
\usepackage{multirow}
\usepackage{hyperref}  
                 \hypersetup{ pdfborder={0 0 0}, 
                              colorlinks=true, 
                              linktoc=page,
                              pdfauthor={Giovanni Staglian{\`o}}, 
                              pdftitle={Some new rational Gushel fourfolds}} 
\usepackage{verbatim}

\theoremstyle{plain} 
\newtheorem{proposition}{Proposition}[section]

\newtheorem*{conjecture}{Conjecture}
\theoremstyle{definition}

\theoremstyle{remark} 
\newtheorem{remark}[proposition]{Remark} 
                                    
\newcommand{\PP}{{\mathbb{P}}}  

\numberwithin{equation}{section}

\title{Some new rational Gushel fourfolds}

\address{Dipartimento di Matematica e Informatica, Universit\`a degli Studi di Catania} 
\author[G. Staglian\`o]{Giovanni Staglian\`o}
\email{giovannistagliano@gmail.com}
\subjclass[2010]{14J35, 
                 14J45, 
                 68W30, 
                 14Q10  
                } 

\begin{document}

\begin{abstract} 
We provide explicit equations and parameterizations of some new rational Gushel-Mukai fourfolds of special type.
\end{abstract}

\maketitle

\section*{Introduction}
The problem of rationality of Fano fourfolds, 
with particular regard to the cases of 
\emph{cubic fourfolds} and \emph{Gushel-Mukai fourfolds},
dates back to classical works; see \emph{e.g.}
\cite{Morin,Fano,Roth1949} (see also Table~\ref{tabella2} below).
Despite the great attention received and the numerous results obtained, 
this is still an open problem for the general fourfold.
For instance, no examples of cubic fourfolds as well as of Gushel-Mukai fourfolds are known to be irrational,
and there are not so many constructions of rational examples.
The contribution of this paper is in the construction of 
some new special examples of rational Gushel-Mukai fourfolds. This is 
achieved by determining their equations through calculations 
with \emph{Macaulay2} \cite{macaulay2}, using mainly the 
packages \emph{SpecialFanoFourfolds} \cite{SpecialFanoFourfoldsSource} and \emph{Cremona} \cite{packageCremona}.
In particular, these packages provide the tools needed to verify the claims in the paper.
\medskip 

Recall that by a result of Mukai \cite{mukai-biregularclassification} (extending 
to all dimensions a result of Gushel proved in \cite{Gu} 
only in dimension three), 
a complex smooth prime Fano fourfold $X$ of degree $10$ and genus $6$, also known as \emph{Gushel-Mukai fourfold},
can be embedded in $\PP^8$
as a quadratic section of a $5$-dimensional linear section $Y\subset \PP^8$ of the 
cone $C(\mathbb{G}(1,4))\subset\PP^{10}$ over the Grassmannian $\mathbb{G}(1,4)\subset\PP^9$ 
of lines in $\PP^4$.
There are two cases:
\begin{itemize}
 \item either $Y$ does not contain the vertex of the cone $C(\mathbb{G}(1,4))$,
 in which case $Y$ is isomorphic to a hyperplane section of $\mathbb{G}(1,4)$, and 
 we have an embedding $\gamma_{X}:X\hookrightarrow \mathbb{G}(1,4)$;
 \item or otherwise $Y$ is isomorphic 
          to a cone over a $4$-dimensional linear section  $Y_0$ of $\mathbb{G}(1,4)$,
          and we have a double cover $\gamma_{X}:X\to Y_0\subset\mathbb{G}(1,4)$.
\end{itemize}
The fourfolds in the second case are called \emph{Gushel fourfolds}, and they are specializations of fourfolds in the first case, called \emph{Mukai (or ordinary) fourfolds}. 
In both cases, the map $\gamma_{X}$ from $X$ to $\mathbb{G}(1,4)$ is called the \emph{Gushel map}.
\medskip 

By results proved in \cite{DIM} (see also \cite{DK1,DK3,DK2}),
Fano fourfolds as above 
are parameterized (up to isomorphism) by the points of a coarse moduli space $\mathcal{M}_{4}$ of dimension $24$,
where the Gushel fourfolds correspond to the points of a closed irreducible subvariety 
$\mathcal{M}_{4}^{G}\subset \mathcal{M}_{4}$ of codimension $2$.
A 
fourfold $[X]\in \mathcal{M}_4$ is said to be \emph{special} (or \emph{Hodge-special}) if 
it 
contains a surface whose cohomology class does not lie in
 $\gamma_{X}^{\ast}(H^4(\mathbb{G}(1,4),\mathbb{Z}))$;
equivalently, $[X]$ is special if and only if $\mathrm{rk}(H^{2,2}(X)\cap H^4(X,\mathbb{Z}))\geq 3$.
A fourfold which corresponds to a very general point in $\mathcal M_4$ or in $\mathcal M_{4}^{G}$ is not special.
More precisely, 
special fourfolds are parametrized by an infinite countable union 
of hypersurfaces $\bigcup_d (\mathcal M_{4})_d\subset \mathcal M_4$, 
labelled by the integers $d\geq 10$ with $d\equiv 0,2$, or $4$ (mod~$8)$.
If $d\equiv 0$ (mod~$4$) then the hypersurface $(\mathcal M_{4})_d$ is irreducible,
while if $d\equiv 2$ (mod~$8$) then it is the union of two irreducible 
components $(\mathcal M_{4})_d'\cup (\mathcal M_{4})_d^{''}$.
When a fourfold 
$[X]$ corresponds to a very general point of a component 
of $(\mathcal M_{4})_d$, then the lattice $ H^{2,2}(X)\cap H^{4}(X,\mathbb Z)$ 
has rank $3$ and discriminant $d$. One says that $X$ has discriminant $d$ if $[X]\in (\mathcal M_{4})_d$.
\medskip 

Suppose we have a fourfold $[X]\in \mathcal M_4$ containing an irreducible 
surface $S$ of degree $\deg(S)$ and (sectional) genus $g(S)$, which has smooth normalization and  only a finite number $\delta$ 
of nodes as singularities. Let 
 $a\sigma_{3,1} + b\sigma_{2,2}$ be the class of $\gamma_{X\ast}(S)$ in the Chow ring of $\mathbb{G}(1,4)$.
 The double point formula (see \cite[Theorem~9.3]{fulton-intersection} and also \cite[Section~7]{DIM}) gives 
 the value of the self-intersection of $S$ in $X$:
\begin{equation}\label{selfIntersection}
(S)_X^2 = 3\,a+4\,b-2\,\deg(S)+4\,g(S)-12\,\chi(\mathcal O_S) + 2\,K_S^2 - 4 + 2\,\delta .
\end{equation}
Thus we have that $[X]\in (\mathcal M_4)_d$, where $d$ is 
the discriminant of the lattice spanned by $(\gamma_X^\ast(\sigma_{1,1}), \gamma_X^{\ast}(\sigma_{2}), [S])$, that is
\begin{equation}\label{discriminant}
d = \mathrm{disc}\left(
     \begin{array}{c|c|c|c}
    & \gamma_X^\ast(\sigma_{1,1}) & \gamma_X^{\ast}(\sigma_{2}) & [S] \\ \hline 
     \gamma_X^\ast(\sigma_{1,1}) & 2 & 2 & b \\ \hline 
      \gamma_X^{\ast}(\sigma_{2}) & 2 & 4 & a \\ \hline 
     \left[S\right] & b & a & (S)_X^2 \\
     \end{array}
     \right)= 4 (S)_X^2 -2 a^{2}+4 a b-4b^{2} .
\end{equation}
Moreover, when $d\equiv 2 \ (\mathrm{mod}\ 8)$, we have that 
$[X]\in(\mathcal M_4)_d'$ if $a+b$ is even, and $[X]\in(\mathcal M_4)_d^{''}$ 
if $b$ is even.

\medskip 

For some values of the discriminant $d$, 
a fourfold $[X]\in (\mathcal M_4)_d$ 
has an \emph{associated} K3 surface of degree $d$;
and for others, it has an \emph{associated} cubic fourfold of discriminant $d$;
see \cite[Section~6]{DIM} for precise definitions.
The first values for which there is an associated K3 surface are:
\begin{equation*}
 { 10},\ {\color{Gray} 12,}\ {\color{Gray} 16,}\ {\color{Gray} 18,}\ { 20},\ {\color{Gray} 24,}\ { 26},\ {\color{Gray} 28,}\ {\color{Gray} 32,}\ { 34},\ {\color{Gray} 36,}\ {\color{Gray} 40,}\ {\color{Gray} 42,}\ {\color{Gray} 44,}\ {\color{Gray} 48,}\ { 50},\ { 52},\ {\color{Gray} 56,}\ { 58},\ {\color{Gray} 60,}\ldots 
\end{equation*}
while, the first values for which there is an associated cubic fourfold are:
\begin{equation*}
 {\color{Gray} 10,}\ { 12},\ {\color{Gray} 16,}\ {\color{Gray} 18,}\ {\color{Gray} 20,}\ {\color{Gray} 24,}\ { 26},\ {\color{Gray} 28,}\ {\color{Gray} 32,}\ {\color{Gray} 34,}\ {\color{Gray} 36,}\ {\color{Gray} 40,}\ {\color{Gray} 42,}\ { 44},\ {\color{Gray} 48,}\ {\color{Gray} 50,}\ {\color{Gray} 52,}\ {\color{Gray} 56,}\ {\color{Gray} 58,}\ {\color{Gray} 60,}\ldots 
\end{equation*}
The notion of associated K3 surface leads to the following conjecture,
which is analogous to the so-called
 \emph{Kuznetsov conjecture} for the rationality of cubic fourfolds (see \cite{kuz4fold,AT,kuz2,Levico,BRS,RS1,RS3}):
\begin{conjecture}
 A fourfold $[X]\in \mathcal M_4$ is rational if and only if 
 it has an associated K3 surface, that is $[X]$ belongs to the infinite union:
 \begin{equation}\label{components}
  (\mathcal M_4)_{10}^{'} \cup (\mathcal M_4)_{10}^{''} \cup (\mathcal M_4)_{20} \cup (\mathcal M_4)_{26}^{'} \cup (\mathcal M_4)_{26}^{''} \cup (\mathcal M_4)_{34}^{'} \cup (\mathcal M_4)_{34}^{''} \cup \cdots 
 \end{equation}
\end{conjecture}
The rationality for fourfolds in $(\mathcal{M}_4)_{10}^{'}$ 
is easy to show (see \cite[Proposition~7.3]{DIM}, \cite[Section~4.4]{ProcRS}, and also \cite[Example~1.1]{HoffSta}), and moreover
the associated K3 surface of degree $10$ appears naturally in the construction 
of the birational map $\PP^4\dashrightarrow X$. The rationality 
for the fourfolds in  $(\mathcal{M}_4)_{10}^{''}$ is classical:
it is achieved by Roth in \cite{Roth1949} as a consequence of a result of Enriques \cite{Enr,EinSh}; 
see also \cite[Proposition~7.5]{DIM}. In this last case, however, the role of the associated K3 surface
is not so clear. 
      In the recent paper \cite{HoffSta}, it is showed that 
a general fourfold in $(\mathcal{M}_4)_{20}$ (and hence \emph{every} by the main result in \cite{KontsevichTschinkelInventiones})
is rational. Again, the associated K3 surface of degree $20$ 
appears in the 
explicit construction leading to rationality.
In conclusion, we have that every fourfold in the first three components of \eqref{components} is rational.

Restricting attention to the case of
 Gushel fourfolds, we point out that 
using the same method presented in \cite{famGushelMukai}, 
just by replacing the role of the 
smooth cubic scroll surface in \cite[Table~2]{famGushelMukai} with that of a cone over a twisted cubic curve,
one can find explicit  Gushel fourfolds in $(\mathcal{M}_4)_{10}^{'}$, $(\mathcal{M}_4)_{10}^{''}$, 
and $(\mathcal{M}_4)_{20}$. In particular, the following three intersections 
\begin{equation*}
  (\mathcal{M}_4)_{10}^{'}\cap \mathcal{M}_4^{G} ,\quad (\mathcal{M}_4)_{10}^{''}\cap \mathcal{M}_4^{G} ,\quad (\mathcal{M}_4)_{20}\cap \mathcal{M}_4^{G} 
\end{equation*}
are not empty, and hence they parametrize  rational fourfolds.

As far as the author knows, no other fourfolds in $ \mathcal M_4$ as well as in $ \mathcal M_4^{G}$ are 
known to be rational.
In the following of this paper, we explain how to find explicit equations and parameterizations 
of rational fourfolds 
in $(\mathcal{M}_4)_{26}''\cap \mathcal M_{4}^{G}$ and
$(\mathcal{M}_4)_{26}''\setminus \mathcal M_{4}^{G}$; see also Table~\ref{tabella1} for a summary.

\subsection*{Acknowledgements}
The author has benefited from discussions with Michele \-{Bolognesi}, Olivier Debarre, and Francesco Russo.

\section{Construction of new rational Gushel fourfolds}\label{sec Gushel}
In this section, we construct rational fourfolds 
in $(\mathcal{M}_4)_{26}''\cap \mathcal M_{4}^{G}$. Here we briefly summarize the construction.
In Subsection~\ref{C26 rat},
starting with a general cubic fourfold $C$ in $\mathcal C_{26}$, 
we explain how to determine the equations of a smooth surface 
$S\subset \PP^8$ of degree $17$ and 
sectional genus $11$, which is isomorphic 
to a triple projection 
of a minimal K3 surface of degree $26$ in $\PP^{14}$. 
In Subsection~\ref{embedding in cone},
we embed the surface $S$ into a smooth quadratic section $X$
of a cone $Y$ over a smooth
$4$-dimensional linear section of $\mathbb{G}(1,4)$.
Then we deduce that $X$ is a fourfold in $(\mathcal{M}_4)_{26}''\cap \mathcal M_{4}^{G}$
since $\gamma_{X \ast}([S]) = 11 \sigma_{3,1} + 6 \sigma_{2,2}$
in the Chow ring of $\mathbb{G}(1,4)$.
In Subsection~\ref{subsec: count parameters 1}, 
we remark that the fourfold $X$ 
may not be a general 
 point in $(\mathcal{M}_4)_{26}''\cap \mathcal M_{4}^{G}$.
In Subsection~\ref{congruence quintics}, 
we illustrate two different methods to deduce that the surface $S$ admits inside the cone $Y$ 
a \emph{congruence of $9$-secant rational normal quintic curves}: through 
a general point of $Y$ there passes a unique rational normal quintic curve 
which is $9$-secant to $S$ and contained in $Y$. 
In Subsection~\ref{particular singular}, we show how from this congruence of quintic curves,
using in an essential way that $Y$ is not smooth,
we can deduce that $X$ is rational.
In Subsection~\ref{summary constr}, we  
describe an explicit birational map between $X$ and the cubic fourfold $C$.

\subsection{Rationality of cubic fourfolds in \texorpdfstring{$\mathcal C_{26}$}{C26}}\label{C26 rat}
Here we recall some results from \cite{RS1,RS3} (see also \cite{ProcRS}) about the rationality of 
special cubic fourfolds of discriminant $26$; 
see also \cite{Hassett,Has00,Levico}  for general facts on cubic fourfolds.

Let $D\subset \PP^5$ be a septimic surface with one node, which is the projection 
of a smooth del Pezzo surface of degree seven $D'\subset\PP^7$ from a general line 
intersecting the secant variety of $D'$ at one point. 
A cubic fourfold $C\subset\PP^5$ containing the surface $D$ 
has discriminant $26$, and more precisely
the locus $\mathcal C_{26}$ of cubic fourfolds of discriminant $26$ 
can be described as the closure inside the moduli space $\mathcal C$ of cubic fourfolds of the locus of fourfolds 
containing such a surface.

The surface $D\subset\PP^5$ admits a \emph{congruence of $5$-secant conics}:
through a general point in $ \PP^5$ there passes a unique $5$-secant conic to $D$.
Moreover, the linear system $|H^0(\mathcal I^2_{D,\PP^5}(5))|$ 
of hypersurfaces of degree $5$ with points of multiplicity $2$ along $D$ 
gives a dominant map 
\begin{equation}\label{KollarMap1}
  \PP^5\dashrightarrow Y_0\subset\PP^7,
\end{equation}
whose general fibers are the conic curves
of the congruence, and where $Y_0$ is a smooth $4$-dimensional linear section of $\mathbb{G}(1,4)\subset\PP^9$.
The restriction of the map \eqref{KollarMap1} to a general cubic fourfold $C$ through $D$ induces 
a birational map $C\dashrightarrow Y_0$, 
whose inverse 
is defined by the linear system $|H^0(\mathcal I^2_{S_0,Y_0}(5))|$ 
of hypersurfaces in $Y_0$ of degree $5$ having points of 
multiplicity $2$ along an irreducible surface $S_0\subset Y_0\subset\PP^7$ 
of degree $17$ and sectional genus $11$ cut out by the $5$ quadrics defining $Y_0$
and $13$ cubics. 

It turns out that $S_0$ is the 
projection of a surface $S\subset\PP^8$
from a special point $p$ on the secant variety of $S$,
where $S$ is a smooth surface of degree $17$,
sectional genus $11$, cut out by $12$ quadrics, and isomorphic 
to a triple projection 
of a minimal K3 surface of degree $26$ in $\PP^{14}$. The equations of $S$ 
can be determined from those of $S_0$
using 
the 
package \emph{IntegralClosure} \cite{IntegralClosureSource};
see also the function \texttt{associatedK3surface} from the package 
\emph{SpecialFanoFourfolds} \cite{SpecialFanoFourfoldsSource}, which does most of this automatically.

\subsection{Gushel fourfolds in \texorpdfstring{$(\mathcal{M}_4)_{26}''$}{(M4)26''}: construction of a triple \texorpdfstring{$S\subset X \subset Y$}{S,X,Y}}\label{embedding in cone}

Continuing from the previous subsection,
let $\nu_p:\PP^8\dashrightarrow\PP^7$ denote the projection from the point $p$ such that $\nu_p(S) = S_0$, and let
$Y=\overline{\nu_p^{-1}(Y_0)}$ be the cone over $Y_0$ of vertex $p$.
Let $X\subset Y$ be a general quadratic section of $Y$ containing $S$.
Then $X$ is a smooth Gushel fourfold which belongs to $(\mathcal M_4)_{26}''$. 

Indeed, the surface $S_0\subset Y_0\subset\mathbb{G}(1,4)$, which is equal to
the image of $S$ via the Gushel map of $X$, has class $11 \sigma_{3,1} + 6 \sigma_{2,2}$
in the Chow ring of $\mathbb{G}(1,4)$, 
as one can  verify by simple calculations with Schubert cycles.
Then, from \eqref{selfIntersection} and \eqref{discriminant} it follows that  
 $(S)_X^2 = 37$ and  $[X]\in  (\mathcal M_4)_{26}^{''}$,
 since one has $\deg(S)=17$, $g(S)=11$, $\chi(\mathcal O_S) = 2$, and $K_S^2 = -1$.
 This calculation can be performed automatically using the functions \texttt{discriminant}
 and \texttt{describe}
from the package \emph{SpecialFanoFourfolds}.
 
\subsection{Count of parameters from the triple \texorpdfstring{$S\subset X\subset Y$}{S,X,Y}}\label{subsec: count parameters 1} 
Let $N_{S/X}$ and $N_{S/Y}$ denote, respectively, 
the normal bundle of the surface $S$ in $X$, and of $S$ in $Y$.
A \emph{Macaulay2} calculation tells us that 
$h^1(N_{S/Y})=0$, $h^0(N_{S/Y})=37$, and  $h^0(N_{S/X})=6$.
It follows that there exists a unique irreducible component $\mathcal S$
of the Hilbert scheme of $Y$ which contains $[S]$, and $\mathcal S$ 
is smooth at $[S]$ of dimension $37$. 
Since we have $h^0(\mathcal I_{S,Y}(2))= 7$ (and this value is minimal on $\mathcal S$),
we deduce by the same
 semicontinuity argument 
explained in \cite[Subsection~1.5]{famGushelMukai} that
inside the $39$-dimensional projective space $\mathbb{P}(H^0(\mathcal O_{Y}(2)))$
of quadratic sections of $Y$ the family of fourfolds containing a surface in $\mathcal S$
has codimension at most $39 - (37 + (7-1) - 6) = 2$. 
This calculation can be performed automatically using the function \texttt{parameterCount} 
from the package \emph{SpecialFanoFourfolds}.

\subsection{Congruence of \texorpdfstring{$9$}{9}-secant quintic curves to \texorpdfstring{$S\subset Y$}{S in Y}}\label{congruence quintics} 
We claim that  the surface $S$ admits inside $Y$ 
a \emph{congruence of $9$-secant rational normal quintic curves}: through 
a general point of $Y$ there passes a unique rational normal quintic curve 
which is $9$-secant to $S$ and contained in $Y$. 

This can be verified by considering the rational map 
$\phi:Y\dashrightarrow \PP^6$ defined by the linear system $|H^0(\mathcal I_{S,Y}(2))|$ 
of quadratic sections through $S$, which is birational onto a non-normal sextic hypersurface $Z\subset \PP^6$.
If $p\in Y$ is a general point, then through the point $\phi(p)\in  Z$ there pass $72$ lines that are contained in $Z$.
These $72$ lines come from  $(2e-1)$-secant curves  to $S$ of degree $e\geq 1$ which pass through $p$ and are contained in $Y$. 
Denoting by $n_e$ the number of such degree-$e$ curves, we have
$n_1=11$, $n_2 = 22$, $n_3 = 32$, $n_4 = 6$, $n_5 = 1$, and $n_e = 0 $ for $e\geq 6$.
This calculation can be performed automatically using the function \texttt{detectCongruence} 
from the package \emph{SpecialFanoFourfolds}.

Alternately, using tools from the package \emph{Cremona}, one verifies that 
 the linear system 
$|H^0(\mathcal I^5_{S,Y}(9))|$ of hypersurfaces in $Y$ 
of degree $9$ with points of multiplicity $5$ along $S$ 
gives a  dominant rational map 
\begin{equation}\label{kollarMap}
 Y\dashrightarrow X'\subset\PP^8,
\end{equation}
whose general fibers 
are rational normal quintic curves, the curves of the congruence to $S$.
The image $X'$ is a Gushel fourfold, a smooth quadratic section of a 
cone $Y'\subset\PP^8$ over a smooth $4$-dimensional linear section  of $\mathbb{G}(1,4)$.

The restriction of the map \eqref{kollarMap}  to $X$ induces a birational map 
\begin{equation}\label{Map}
Y \supset X \stackrel{\simeq}{\dashrightarrow}  X' \subset Y' ,
\end{equation}
whose inverse is of the same type, \emph{i.e}, it is
the restriction to $X'$ of the 
rational map $Y'\dashrightarrow X$ defined 
by the linear system of 
hypersurfaces in $Y'$
of degree $9$ with points of multiplicity $5$ along a smooth triple projection $S'\subset Y' $ 
of a minimal K3 surface of degree $26$. 
In particular, the Gushel fourfold $X'$ also belongs to $(\mathcal M_4)_{26}^{''}$.

\subsection{Rationality of \texorpdfstring{$X$}{X} from the rationality of a particular rational singular fourfold \texorpdfstring{$\hat{X}$}{Xh} with \texorpdfstring{$S\subset \hat{X}\subset Y$}{S in X in Y}}\label{particular singular}
Let $\hat{X}\subset Y$ be a general quadratic section of $Y$
 containing the surface $S$ and the vertex $p$ of $Y$. Then $\hat{X}$ has $p$
as the only singularity and the restriction of the projection from $p$ 
induces a birational map $\hat{X}{\dashrightarrow} Y_0$, whose inverse is defined by the quadrics
through a minimal K3 surface $F$ of degree $10$. 
Since $Y_0=\mathbb{G}(1,4)\cap \PP^7$ is rational (indeed, the projection from the unique $\sigma_{2,2}$-plane contained in it gives a birational map onto $\PP^4$), we have that  $\hat{X}$ is rational.
On the other hand, the restriction 
of the map \eqref{kollarMap}  induces a birational map 
\[
 Y \supset \hat{X} \stackrel{\simeq}{\dashrightarrow} X' \subset Y',
\]
due to the fact that $\hat{X}$ is \emph{transversal} to the congruence to $S$: 
the quintic curve of the congruence passing through a general point of $\hat{X}$ 
is not contained in $\hat{X}$.
Therefore we deduce that also $X'$ and hence $X$ are rational. 

\subsection{Summary construction}\label{summary constr}
Summing up we have the following diagram of birational maps, connecting explicitly
general cubic fourfolds in $\mathcal C_{26}$  
to Gushel fourfolds in $(\mathcal M_{4})_{26}''$:
 \begin{equation}\label{mappe}
 \xymatrix{
 {\color{Gray} \PP^5} \ar@/^1.9pc/@{-->}@[Gray][rrd] & & & & {\color{Gray} Y}\ar@/^1.9pc/@{-->}@[Gray][rrd] & & {\color{Gray} Y'}\ar@/^1.9pc/@{-->}@[Gray][rrd] \\
  C\ar@{^{(}->}@[Gray][u]\ar@/^1.0pc/@{-->}[rr]^{|H^0(\mathcal{I}^{2}_{D}(5))|} && Y_0\ar@/^1.0pc/@{-->}[ll]^{|H^0(\mathcal{I}^{2}_{S_0}(5))|} 
  \ar@/^1.0pc/@{-->}[rr]^{|H^0(\mathcal{I}_{F}(2))|} && \hat{X}\ar@{^{(}->}@[Gray][d]\ar@{^{(}->}@[Gray][u]\ar@/^1.0pc/@{-->}[ll]^{|H^0(\mathcal{I}_{p}(1))|} 
  \ar@/^1.0pc/@{-->}[rr]^{|H^0(\mathcal{I}^{5}_{S}(9))|} && X'\ar@{^{(}->}@[Gray][u]\ar@{^{(}->}@[Gray][d]\ar@/^1.0pc/@{-->}[ll]^{|H^0(\mathcal{I}^{5}_{\hat{S}'}(9))|} 
  \ar@/^1.0pc/@{-->}[rr]^{|H^0(\mathcal{I}^{5}_{S'}(9))|} && X\ar@{^{(}->}@[Gray][d]\ar@/^1.0pc/@{-->}[ll]^{|H^0(\mathcal{I}^{5}_{S}(9))|} \\
  & & & & { \color{Gray} Y } \ar@/^1.9pc/@{-->}@[Gray][llu] & & {\color{Gray} Y' } \ar@/^1.9pc/@{-->}@[Gray][llu] & & { \color{Gray} Y} \ar@/^1.9pc/@{-->}@[Gray][llu]
 } 
\end{equation}
where $[C]\in \mathcal C_{26}$ 
is a general cubic fourfold of discriminant $26$ containing 
a septimic one-nodal del Pezzo surface $D$;
$[X]$, $[X']\in (\mathcal M_{4})_{26}''\cap \mathcal M_{4}^{G}$ are (smooth) Gushel fourfolds 
of discriminant $26$ contained, respectively, in cones $Y$ and $Y'$ over $Y_0=\mathbb{G}(1,4)\cap \PP^7$; 
$\hat{X}\subset Y$ is a Gushel fourfold singular at $p$, where 
$p$ is the vertex of $Y$;
$S\subset X\cap \hat{X}$ and $S',\hat{S}'\subset X'$
are smooth surfaces, isomorphic to triple projections of minimal K3 surfaces of degree $26$ (the intersection $S'\cap \hat{S}'$ 
consists of a twisted cubic curve and $31$ points);
$S_0\subset Y_0$ is the projection of $S$ from $p$;
and $F\subset Y_0$ is a minimal K3 surface
of degree $10$.

\subsection{Ancillary files}\label{computation 1}
We provide an ancillary file, named \texttt{gushel26.m2},
containing explicit equations for an example of map as \eqref{Map} over the finite field $\mathbb{F}_{10000019}$ 
(this is only necessary to reduce the size of the file).
After loading that file in \emph{Macaulay2}, 
some variables will be defined as follows:
\begin{description}
 \item[\texttt{X, X'}] two instances of the type \texttt{SpecialGushelMukaiFourfold}, respectively, the source and target of the map \eqref{Map};
 \item[\texttt{psi, psi'}] two instances of the type \texttt{RationalMap}, respectively,
 the birational map \eqref{Map} and its inverse (\texttt{psi'} is the same that \texttt{inverse psi});
 \item[\texttt{Psi, Psi'}] two instance of the type \texttt{RationalMap}, respectively,
 the dominant rational map \eqref{kollarMap} from $Y$ to $X'$ that extends \texttt{psi},
 and the analogous map from $Y'$ to $X$ that extends \texttt{psi'}. 
\end{description}
For technical details about these types of data, we refer to the documentation of the packages \emph{Cremona} \cite{packageCremona} and 
\emph{SpecialFanoFourfolds} \cite{SpecialFanoFourfoldsSource}.
We now show how to load the file 
and extract some basic information from it (some output lines are omitted for brevity).
\begin{tcolorbox}[breakable=true,boxrule=0pt,opacityback=0.0,enhanced jigsaw]
{\footnotesize
\begin{Verbatim}[commandchars=&\[\]]
&colore[darkorange][$ M2 --no-preload]
&colore[output][Macaulay2, version 1.16]
&colore[darkorange][i1 :] &colore[airforceblue][needs] "&colore[bleudefrance][gushel26.m2]"
&colore[darkorange][i2 :] &colore[airforceblue][describe] X
&colore[circOut][o2 =] &colore[output][Special Gushel-Mukai fourfold of discriminant 26('')]
&colore[output][     containing a surface in PP^8 of degree 17 and sectional genus 11]
&colore[output][     cut out by 12 hypersurfaces of degree 2]
&colore[output][     and with class in G(1,4) given by 11*s_(3,1)+6*s_(2,2)]
&colore[output][     Type: Gushel (not ordinary)]
&colore[darkorange][i3 :] &colore[airforceblue][describe] Psi
&colore[circOut][o3 =] &colore[output][rational map defined by forms of degree 9]
&colore[output][     source variety: 5-dimensional variety of degree 5 in PP^8]
&colore[output][                     cut out by 5 hypersurfaces of degree 2]
&colore[output][     target variety: 4-dimensional variety of degree 10 in PP^8]
&colore[output][                     cut out by 6 hypersurfaces of degree 2]
&colore[output][     dominance: true]
&colore[output][     projective degrees: {5, 45, 211, 200, 50, 0}]
&colore[darkorange][i4 :] S = &colore[airforceblue][first ideals] X; &colore[commentoutput][-- the surface S]
&colore[darkorange][i5 :] h = &colore[airforceblue][detectCongruence](X,5); &colore[commentoutput][-- the congruence to S]
&colore[darkorange][i6 :] p = &colore[airforceblue][point] X; &colore[commentoutput][-- random point on X]
&colore[darkorange][i7 :] C = h p; &colore[commentoutput][-- 9-secant quintic curve to S passing through p]
\end{Verbatim}
} 
\end{tcolorbox}

\section{Rationality 
without passing through singular fourfolds}\label{particular smooth}
A count of parameters shows that 
a general triple projection 
 of a general minimal K3 surface of degree $26$ 
is contained in a one-dimensional family 
of $5$-dimensional linear sections  of cones $C(\mathbb{G}(1,4))\subset\PP^{10}$.
So one expects that the Gushel fourfolds constructed in Section~\ref{sec Gushel} 
can be deformed to ordinary fourfolds. Anyway, the simple argument 
given in Subsection~\ref{particular singular} to deduce the rationality of the fourfold $X$
is not available in the ordinary case. 
In this section, we remedy this with a construction that does not involve singular Gushel fourfolds.

Indeed, keeping the notation as in Section~\ref{sec Gushel}, here we construct  
a special Gushel fourfold $\tilde{X}$ with $S\subset \tilde{X}\subset Y$,
uniquely determined by the embedding $S\subset Y$, which turns out to 
be smooth and transversal to the congruence of $9$-secant quintic curves to $S$.
We then show that $\tilde{X}$ contains another ``simpler'' surface $T$ which admits 
inside $Y$ a congruence of $5$-secant  cubic curves. So we deduce that $\tilde{X}$ is rational,
and hence that  $X$ is rational.

\subsection{Construction of \texorpdfstring{$\tilde{X}$}{X} and of a dominant map \texorpdfstring{$\tilde{X}\dashrightarrow\PP^2$}{X-->P2}}
Consider again the birational map introduced in Subsection~\ref{congruence quintics},
\[
\phi:Y\dashrightarrow Z\subset \PP^6 ,
\]
defined by the linear system $|H^0(\mathcal I_{S,Y}(2))|$, 
where $Z=\overline{\phi(Y)}\subset\PP^6$ is a sextic hypersurface.
We take $\tilde{X}$ to be the top dimensional component of the (closure of the) exceptional locus of $\phi$, and the map $\tilde{X}\dashrightarrow\PP^2$ to be defined by the quadrics through the (closure of the) union of all $4$-secant conics to $S$ contained in $\tilde{X}$.
Let us provide some more detail.
\subsubsection{The fourfold $\tilde{X}$}
Using tools from the 
package \emph{Cremona}, one verifies that the base locus scheme of the inverse map $\phi^{-1}:Z\dashrightarrow Y$ 
is the union of the following components:
\begin{itemize}
 \item a smooth cubic fourfold $\tilde{C}\subset\PP^5\subset\PP^6$ (with $[\tilde{C}]\in \mathcal{C}_{26}$, see Remark~\ref{remC26}); the fiber of $\phi$ at a general point of $\tilde{C}$ 
 consists of two points;
 \item a smooth cubic scroll surface $\Sigma_3\subset\PP^4\subset\PP^6$, which is double; the fiber of $\phi$ at a general point 
 of $\Sigma_3$ is an irreducible conic curve which is $4$-secant to $S$;
 \item a  surface of degree $46$ cut out in $\PP^6$ by $81$ quintic hypersurfaces; 
 the fiber of $\phi$ at a general point 
 of this surface is a $2$-secant line to $S$.
\end{itemize}
Then, the  fourfold $\tilde{X}$ is taken to be $\overline{\phi^{-1}(\tilde{C})}$, so that 
 the restriction of $\phi$
induces a generically finite map of degree $2$ from $\tilde{X}$ to $\tilde{C}$.
One 
sees
that $\tilde{X}$ is a smooth Gushel fourfold which is transversal to the congruence to $S$ 
and in particular it is birational to $X'$ via the restriction of the map \eqref{kollarMap}.
Moreover, the inverse map $X'\dashrightarrow \tilde{X}$ 
is defined once again by the linear system $|H^0(\mathcal I^{5}_{\tilde{S}'}(9))|$,
where $\tilde{S}'\subset X'\subset Y'$
is a smooth triple projection of a minimal K3 surface of degree $26$. We stress that
this surface $\tilde{S}'\subset Y'$ and the fourfold $\tilde{X}$
  are uniquely determined by the embedding $S\subset Y$.

\subsubsection{The map $\tilde{X}\dashrightarrow\PP^2$}
The intersection $\tilde{C}\cap\Sigma_3$ is a twisted cubic curve, and the fiber of $\phi$ 
at a general point of this curve is a $4$-secant conic curve to $S$ contained in $\tilde{X}$.
The inverse image $R=\overline{\phi^{-1}(\tilde{C}\cap\Sigma_3)}$ is an irreducible surface contained in $\tilde{X}$
of degree $17$ and sectional genus $6$ 
cut out in $\PP^8$ 
by $9$ quadrics and $7$ cubics.
Let 
\begin{equation}\label{MapEta}
     \eta:Y\dashrightarrow \PP^3
\end{equation}
be the rational map 
defined by the linear system $|H^0(\mathcal I_{R,Y}(2))|$ of quadratic sections of $Y$ through $R$. 
Since  $R\subset\tilde{X}$, the restriction of $\eta$ induces 
another map
\begin{equation}\label{MapEta2}
     \eta|_{\tilde{X}}:\tilde{X}\dashrightarrow \PP^2 \subset\PP^3 .
\end{equation}
One sees that $\eta$ (resp., $\eta|_{\tilde{X}}$) is a dominant
map, 
whose general fibers are surfaces like $S$, that is, 
smooth triple projections 
of minimal K3 surfaces of degree $26$.
Moreover the surface $S$ 
is recovered as a special fiber, and
all these fibers share the same twisted cubic curve contained in $S$.

\begin{remark}
 The projection of the surface $S$ from the plane spanned by a general $4$-secant
 conic to $S$ contained in $Y$ (computable as the the 
 fiber of $\phi$ at a general point on $\Sigma_3$) is a smooth surface in $\PP^5$ 
 of degree $13$ and sectional genus $11$ cut out by $6$ cubics.
 This surface in $\PP^5$ 
 admits a congruence of $14$-secant rational normal quintic curves 
 from which one can deduce the rationality 
 for cubic fourfolds of discriminant $26$.
 (See also the example in \cite[Table~1, row~15]{RS3} of a nodal surface in $\PP^5$.)
\end{remark}

\subsection{Special fibers of the map \texorpdfstring{$\tilde{X}\dashrightarrow \PP^2$}{X-->P2}: construction of a surface \texorpdfstring{$T\subset\tilde{X}\subset Y$}{T in X in Y}}\label{Gushel fourfold T,Xtilde}
The $3$-dimensional projective space, image of the map $\eta$, contains a special plane $\Pi$ 
which intersects $\overline{\eta(\tilde{X})}\simeq\PP^2$ along a line $L$.
The fiber  of $\eta$ at a general point of $L$ (as well as of $\Pi$)  
is an irreducible rational surface $T\subset Y\subset \PP^8$ of degree $11$, sectional genus $3$,
cut out in $\PP^8$ by $16$ quadrics, having smooth normalization and 
 a node as the only singularity. 
 In the Chow ring of $\mathbb{G}(1,4)$ we have 
 $[\gamma_{\tilde{X}}(T)] = 7 \sigma_{3,1} + 4 \sigma_{2,2}$, so
 from  \eqref{selfIntersection} it follows that $(T)_{\tilde{X}}^2 = 19$ (this is also confirmed from the fact that two 
 fibers of $\eta$ corresponding to two general points of $L$ intersects at $19$ points), 
 and \eqref{discriminant} tells us that
 any smooth quadratic section of $Y$ containing $T$ is a Gushel fourfold 
 of discriminant $26$, hence corresponding to a point of~$(\mathcal M_4)_{26}^{''}$.
 
 From another point of view, using the map $\eta$,
 we are able to obtain a degeneration of the surface $S$ as the union of the surface 
 $T$ and a smooth surface $Q$ of degree $6$ and sectional genus $2$ 
 with $[\gamma_{\tilde{X} \ast}(Q)] = 2(2\sigma_{3,1} + \sigma_{2,2}) = 2 \sigma_2 \sigma_1^2$ and such that 
 the intersection $T\cap Q$ is an irreducible curve of degree $7$ with $p_a=2$ and $p_g=1$.
 This surface $Q$ is contained in the base locus of $\eta$. More precisely,
 the support of the base locus of $\eta$ is the union of $Q$ with the surface $R$.

\begin{remark}\label{remC26} 
The plane $\Pi$ 
 can be calculated as the image 
 via $\eta$
of the $3$-dimensional linear space $\overline{\nu_p^{-1}(P)}$, where $\nu_p:Y\dashrightarrow Y_0$
is the projection from the vertex of $Y$, and $P$ is the unique $\sigma_{2,2}$ plane contained in $Y_0$.
Moreover,
the intersection $\overline{\nu_p^{-1}(P)}\cap \tilde{X}$ 
 is a quadric surface which is sent birationally by 
 $\phi$ to a one-nodal septimic surface $\tilde{D}\subset\tilde{C}$ 
 as the surface $D\subset C$ considered in Subsection~\ref{C26 rat}. In particular, we also deduce that $\tilde{C}$
 is a cubic fourfold of discriminant $26$.
\end{remark}

\subsection{Count of parameters from the triple \texorpdfstring{$T\subset \tilde{X}\subset Y$}{T,X,Y}} 
As in Subsection~\ref{subsec: count parameters 1}, we compute with \emph{Macaulay2} that 
$h^1(N_{T/Y})=0$, $h^0(N_{T/Y})=29$, and  $h^0(N_{T/\tilde{X}})=2$.
Assuming that the Hilbert scheme $\mathrm{Hilb}_Y$ of $Y$ is smooth at $[T]$ 
(which is reasonable, but not guaranteed since $T$ is not a local complete intersection), we have that 
$\mathrm{Hilb}_Y$ contains a unique irreducible component $\mathcal T$
 which contains $[T]$, and the dimension of $\mathcal T$ 
is $29$. 
Therefore, since we have $h^0(\mathcal I_{T,Y}(2))= 11$ and this value is minimal,
we deduce that
inside the projective space $\mathbb{P}(H^0(\mathcal O_{Y}(2)))$
of quadratic sections of $Y$ the family of fourfolds containing a surface in $\mathcal T$
has codimension at most $39 - (29 + (11-1) - 2) = 2$.
 
 \subsection{Congruence of \texorpdfstring{$5$}{5}-secant cubic curves to \texorpdfstring{$T\subset Y$}{T in Y}  and rationality of \texorpdfstring{$\tilde{X}$}{X}}
 
 The surface $T\subset Y$  admits inside $Y$
 a congruence of $5$-secant twisted cubic curves,
 and $\tilde{X}$ is transversal to this congruence.

 Indeed one verifies that the linear system $|H^0(\mathcal I_{T,Y}^3(5))|$ 
 of 
hypersurfaces in $Y$
of degree $5$ with points of multiplicity $3$ along 
   $T$ gives a dominant rational map 
 \begin{equation}\label{MapFromYtoW}
  Y\dashrightarrow W\subset\PP^{10}
 \end{equation}
 onto a smooth $4$-dimensional linear section $W$ of $\mathbb{G}(1,5)\subset\PP^{14}$,
 and whose general fibers are twisted cubic curves. The restriction 
 of \eqref{MapFromYtoW} to $\tilde{X}$ induces a birational map 
 \begin{equation}\label{MapFromXtildetoW}
  Y\supset \tilde{X}\stackrel{\simeq}{\dashrightarrow} W\subset\PP^{10} ,
 \end{equation}
 whose inverse is defined by the linear system of 
 hypersurfaces of degree $5$ with points
 of multiplicity $3$ along a smooth surface $U\subset W$ of degree $21$ and sectional genus $13$, isomorphic 
 to a double projection of a simple projection of a minimal K3 surface of degree $26$.
 
 We deduce the rationality of $\tilde{X}$ from that of $W$. Indeed, 
 $W$ must contain a quintic del Pezzo surface, and 
 it is classically known that
the linear system of hyperplanes through this surface gives a birational map $W\dashrightarrow\PP^4$.

 The congruence to $T\subset Y$ can be also verified by considering the 
 map 
 $
  Y\dashrightarrow \PP^{10},
 $
 defined by the linear system $|H^0(\mathcal I_{T,Y}(2))|$ of quadratic sections of $Y$ through $T$, 
 which turns out to be  birational onto a fivefold of degree $20$ cut out by $7$ quadrics.
 Through the general point of this fivefold there pass $12$ lines, which come 
 from seven $1$-secant lines to $T$, four $3$-secant conics to $T$, 
 and one single $5$-secant twisted cubic to $T$.

 There is a further way to find the congruence to $T\subset Y$. Indeed,
 one has that the reducible surface $T\cup Q$ 
 considered in Subsection~\ref{Gushel fourfold T,Xtilde}, which is a degeneration of the surface $S$,
 admits  a congruence of $9$-secant quintic curves, exactly as $S$ does.
 In this degenerate case, the curves of the congruence
  split into $4$-secant conics to $T\cup Q$ and $5$-secant twisted cubics to $T$.  
 
 \subsection{Summary construction} With the notation above introduced, we have the 
 following diagram involving 
 cubic fourfolds in $\mathcal{C}_{26}$ and Gushel fourfolds in $(\mathcal M_4)_{26}^{''}$, and where 
 all the fourfolds are smooth. 
  \begin{equation}\label{mappe2}
 \xymatrix{
  \tilde{C} & & {\color{Gray} Y}\ar@/^1.9pc/@{-->}@[Gray][rrd] & & {\color{Gray} Y'}\ar@/^1.9pc/@{-->}@[Gray][rrd] \\
   W 
  \ar@/^1.0pc/@{-->}[rr]^{|H^0(\mathcal{I}^{3}_{U}(5))|} && \tilde{X}\ar@/_2.0pc/@{-->}[llu]_{2:1}^{H^0(\mathcal{I}_{S}(2))}\ar@{^{(}->}@[Gray][d]\ar@{^{(}->}@[Gray][u]\ar@/^1.0pc/@{-->}[ll]^{|H^0(\mathcal{I}^{3}_{T}(5))|} 
  \ar@/^1.0pc/@{-->}[rr]^{|H^0(\mathcal{I}^{5}_{S}(9))|} && X'\ar@{^{(}->}@[Gray][u]\ar@{^{(}->}@[Gray][d]\ar@/^1.0pc/@{-->}[ll]^{|H^0(\mathcal{I}^{5}_{\tilde{S}'}(9))|} 
  \ar@/^1.0pc/@{-->}[rr]^{|H^0(\mathcal{I}^{5}_{S'}(9))|} && X\ar@{^{(}->}@[Gray][d]\ar@/^1.0pc/@{-->}[ll]^{|H^0(\mathcal{I}^{5}_{S}(9))|} \\
   & & { \color{Gray} Y } \ar@/^1.9pc/@{-->}@[Gray][llu] & & {\color{Gray} Y' } \ar@/^1.9pc/@{-->}@[Gray][llu] & & { \color{Gray} Y} \ar@/^1.9pc/@{-->}@[Gray][llu]
 } 
\end{equation}
 
 \subsection{Ancillary files}\label{computation 2}
 Continuing from Subsection~\ref{computation 1}, 
 after the loading of the ancillary file in \emph{Macaulay2},
 some other variables are defined as following:
 \begin{description}
  \item[\texttt{Xtilde}] an instance of the type \texttt{SpecialGushelMukaiFourfold}, 
  the Gushel fourfold corresponding to the pair $(T,\tilde{X})$ 
  constructed above;
  \item[\texttt{Eta, eta}] two instances of the type \texttt{RationalMap}, respectively, 
  the map $Y\dashrightarrow\PP^3$ in \eqref{MapEta} and its restriction 
  $\tilde{X}\dashrightarrow\PP^2\subset\PP^3$ in \eqref{MapEta2}.
 \item[\texttt{Psi2, psi2}] two instances of the type \texttt{RationalMap}, respectively, 
 the map $Y\dashrightarrow W$ in \eqref{MapFromYtoW} and its restriction $\tilde{X}\dashrightarrow W$ 
 in \eqref{MapFromXtildetoW}.
 \end{description}
By way of example, we now compute a general and a special fiber of 
the map \eqref{MapEta2}.
\begin{tcolorbox}[breakable=true,boxrule=0pt,opacityback=0.0,enhanced jigsaw]
{\footnotesize
\begin{Verbatim}[commandchars=&\[\]]
&colore[darkorange][i8 :] p = &colore[airforceblue][point image] eta; &colore[commentoutput][-- random point on the image of eta]
&colore[darkorange][i9 :] L = &colore[airforceblue][ideal] (&colore[airforceblue][gens target] eta)_{2,3}; &colore[commentoutput][-- the special line L in P^3 (Subsect. 2.2)]
&colore[darkorange][i10 :] q = &colore[airforceblue][point] L; &colore[commentoutput][-- random point on L]
&colore[darkorange][i11 :] &colore[airforceblue][describe](eta||p) &colore[commentoutput][-- fiber at p]
&colore[circOut][o11 =] &colore[output][source variety: surface of degree 17 and sectional genus 11 in PP^8]
&colore[output][                      cut out by 12 hypersurfaces of degree 2]
&colore[output][      target variety: one-point scheme in PP^3]
&colore[darkorange][i12 :] &colore[airforceblue][describe](eta||q) &colore[commentoutput][-- fiber at q]
&colore[circOut][o12 =] &colore[output][source variety: surface of degree 11 and sectional genus 3 in PP^8]
&colore[output][                      cut out by 16 hypersurfaces of degree 2]
&colore[output][      target variety: one-point scheme in PP^3]
\end{Verbatim}
}
\end{tcolorbox}
 
 \section{Construction of new rational Mukai fourfolds}\label{sec Mukai}
 In this section, we briefly explain how using 
 the same method given in \cite{famGushelMukai}, one can
directly  construct a rational one-nodal surface  $T\subset Y$ of degree $11$ and 
sectional genus $3$
as the one constructed in Section~\ref{particular smooth}. Actually, we do better than this:
We are able to construct a 
rational one-nodal surface  $\breve{T}$ of degree $11$ and 
sectional genus $3$
inside a smooth hyperplane section $\breve{Y}$ of $\mathbb{G}(1,4)$.
 This leads us to find out  rational ordinary fourfolds in $(\mathcal M_4)_{26}^{''}$.
 
\subsection{Construction of \texorpdfstring{$\breve{T}\subset\breve{Y}$}{T in Y}}
Let $E'\subset\PP^6$ be the image of the plane via the linear system 
of quartic curves with one double point $q_0$ and $5$ simple base points $q_1,\ldots,q_5$. 
Then $E'$ 
is a smooth surface of degree $7$ and sectional genus $2$ cut out by $8$ quadrics.\footnote{The 
surface $E'\subset\PP^6$ 
is a hyperplane section of a so-called \emph{Edge variety} \cite{Edge}.
This surface also occurs in the classification of \emph{special birational transformations},
see \cite[Table~4, Case~VI]{Sta19}.}
Let $E\subset\PP^5$ be the projection of $E'$ from a general point 
on the secant variety of $E'$. 
Then $E$ is an irreducible surface of degree 
$7$, sectional genus $2$, cut out by $2$ quadrics and $5$ cubics,
with a single node as the only singularity, and having normalization
isomorphic to $E'$. 
Take $B\subset\PP^5$ to be 
a general smooth cubic scroll surface which cuts $E$ 
along a quintic elliptic curve (such a curve is obtained as the image on $E'$ and hence on $E$ of a general plane cubic curve
passing through the six base points $q_0,\ldots,q_5$).

The linear system of quadrics through $B$ defines a birational map 
\begin{equation}\label{SempleMap}
 \PP^5\dashrightarrow \breve{Y}\simeq \mathbb{G}(1,4)\cap\PP^8\subset\PP^8 
\end{equation}
onto a smooth hyperplane section $\breve{Y}$ of $\mathbb{G}(1,4)\subset\PP^9$.
The restriction of this map to $E$ induces an isomorphism 
between $E$ and a surface $\breve{T}\subset\breve{Y}$.

This surface $\breve{T}\subset\breve{Y}$  shares all the properties of the surface $T\subset Y$
constructed in Section~\ref{particular smooth}.
In particular, $\breve{T}\subset\breve{Y}$ is 
a one-nodal irreducible surface 
of degree $11$,
sectional genus $3$,
cut out in $\PP^8$ by $16$ quadrics,
and having class $7 \sigma_{3,1} + 4 \sigma_{2,2}$ 
 in the Chow ring of $\mathbb{G}(1,4)$.
 A general quadratic section of $\breve{Y}$ through $\breve{T}$ 
 gives an ordinary Gushel-Mukai fourfold in $(\mathcal M_4)_{26}^{''}$,
 and such a fourfold is birational to a smooth $4$-dimensional 
 linear section of $\mathbb{G}(1,5)\subset\PP^{14}$
 via the linear system of quintic hypersurfaces 
 with triple points along $\breve{T}$.

\begin{remark}
By counting parameters, one sees that the family of the 
reducible surfaces $B\cup E\subset\PP^5$ has dimension $48$, and hence 
 the family of the surfaces $\breve{T}\subset \breve{Y}$ obtained by this construction 
has dimension $48 - 35 + 15 = 28 < 29 = h^0(N_{\breve{T}/\breve{Y}})$;
see also \cite[Remark~2.4]{famGushelMukai}.
\end{remark}

\begin{remark}
 Our ancillary file (see Subsections~\ref{computation 1} and \ref{computation 2})
 also provides an explicit example 
 of ordinary Gushel-Mukai fourfold corresponding to the pair
  $(\breve{T},\breve{X})$, where $\breve{X}$ is a general quadratic section of $\breve{Y}$ through $\breve{T}$.
 The file includes 
 the parameterization of $E'$, the nodal projection of $E'$ onto $E$, the map \eqref{SempleMap},
 the dominant map $\breve{Y}\dashrightarrow \mathbb{G}(1,5)\cap\PP^{10}$ 
 whose general fibers are $5$-secant twisted cubic curves to $\breve{T}$,
 and the restriction of this map 
 to $\breve{X}$ together with its inverse.
\end{remark}

\section{Summary table of examples}
Table~\ref{tabella1} summarizes some information about the 
   Gushel-Mukai fourfolds known to be rational. 
 It includes all the examples we found in the literature,
 the example constructed in the present paper,
 and some few others of which we omit the details. 
 In all cases, except the 4th, the rationality of the fourfold follows from 
  a congruence  of $(2e-1)$-secant curves of degree $e\geq1$.
\begin{table}[htbp]
\renewcommand{\arraystretch}{1.0} 
\centering
\tabcolsep=0.3pt 
\begin{adjustbox}{width=\textwidth}
\begin{tabular}{|c|c|c|c|c|c|c|c|c|c|}
\hline
\rowcolor{gray!5.0}
Ref. & {\begin{tabular}{c} Surface $S\subset Y$  \end{tabular}} & $K_S^2$ & {\begin{tabular}{c} Class in \\ $\mathbb{G}(1,4)$\end{tabular}} & {\begin{tabular}{c} Locus in $\mathcal M_4$ \end{tabular}}  & $h^0(\mathcal I_{S/Y}(2))$ & $h^0(N_{S/Y})$ & $h^0(N_{S/X})$ &
  \scriptsize{\begin{tabular}{c} Curves of degree $e$ in $Y$\\ passing though a general \\ point of $Y$ and that  are \\ $(2e-1)$-secant to $S$ \\ for $e=1,2,3,4,5$  \end{tabular}} \\
\hline \hline 
\cite{DIM} & {\begin{tabular}{c} Quadric surface   \end{tabular}} & $8$  & $\sigma_{3,1}+\sigma_{2,2}$ & $(\mathcal M_4)_{10}^{'}$ & $31$ & $8$ & $0$ & $1$, $0$, $0$, $0$, $0$  \\
\hline
\cite{ProcRS}  &  {\begin{tabular}{c} K3 surface of degree \\ $14$ and genus $8$  \end{tabular}}  & $0$  & $9 \sigma_{3,1}+ 5 \sigma_{2,2}$ & $(\mathcal M_4)_{10}^{'}$ & $10$ & $39$ & $10$ & $9$, $8$, $1$, $0$, $0$   \\
\hline 
\cite{Roth1949} & {\begin{tabular}{c} Plane   \end{tabular}}  & $9$  & $\sigma_{3,1}$ & \begin{tabular}{c} codim $1$ in \\ $(\mathcal M_4)_{10}^{'}$ \end{tabular} & $34$ & $4$ & $0$ & $1$, $0$, $0$, $0$, $0$ \\
\hline 
\cite{Roth1949} & {\begin{tabular}{c} Quintic del Pezzo \\ surface   \end{tabular}}  & $5$  & $3 \sigma_{3,1}+2 \sigma_{2,2}$ & $(\mathcal M_4)_{10}^{''}$ & $24$ & $18$ & $3$ & $3$, $0$, $0$, $0$, $0$  \\
\hline
\cite{HoffSta} & {\begin{tabular}{c} Rational surface of \\ degree $9$ and genus $2$  \end{tabular}}  & $5$ & $6 \sigma_{3,1}+3 \sigma_{2,2}$ & $(\mathcal M_4)_{20}$ & $14$ & $25$ & $0$ & $6$, $1$, $0$, $0$, $0$ \\
\hline 
  & {\begin{tabular}{c} Rational surface of \\ degree $13$ and genus $6$  \end{tabular}}  & $-2$ & $8 \sigma_{3,1}+5 \sigma_{2,2}$ & \begin{tabular}{c} locus in \\ $(\mathcal M_4)_{20}$ \end{tabular} & $10$ & $33$ & $4$ & $8$, $7$, $1$, $0$, $0$ \\
\hline 
\cite{famGushelMukai} & {\begin{tabular}{c} Septic scroll  \end{tabular}}  & $8$ & $4 \sigma_{3,1} + 3 \sigma_{2,2}$ & \begin{tabular}{c} codim $2$ in \\ $(\mathcal M_4)_{20}$ \end{tabular} & $16$ & $21$ & $0$ & $4$, $1$, $0$, $0$, $0$ \\
\hline 
\S \ref{sec Gushel} & {\begin{tabular}{c} Triple projection of a \\ minimal K3 surface of \\ degree $26$ and genus $14$  \end{tabular}}  & $-1$ & $11 \sigma_{3,1} + 6 \sigma_{2,2}$ & \begin{tabular}{c} locus in \\ $(\mathcal M_4)_{26}^{''}\cap \mathcal M_4^{G}$ \end{tabular} & $7$ & $37$ & $6$ & $11$, $22$, $32$, $6$, $1$ \\
\hline 
\S \ref{particular smooth},\ref{sec Mukai} & {\begin{tabular}{c} Rational $1$-nodal surface \\ of degree $11$ and genus $3$  \end{tabular}}  & $3$ & $7 \sigma_{3,1} + 4 \sigma_{2,2}$ & \begin{tabular}{c} locus in \\ $(\mathcal M_4)_{26}^{''}$ \end{tabular} & $11$ & $29$ & $2$ & $7$, $4$, $1$, $0$, $0$ \\
\hline 
\end{tabular}
\end{adjustbox}
 \caption{Rational Gushel-Mukai fourfolds $X$ 
 obtained as quadratic sections of 
  $Y\simeq C(\mathbb{G}(1,4))\cap\PP^8$ through  surfaces $S\subset Y$.}
 \label{tabella1} 
\end{table}

\begin{table}[htbp]
\centering 
\tabcolsep=1.2pt 
\begin{adjustbox}{width=\textwidth}
\begin{tabular}{cccc}
\hline
\rowcolor{gray!5.0}
Fourfold & \begin{tabular}{c} Irrationality of \\ very general \end{tabular} & \begin{tabular}{c} Description of \\ the rational ones \end{tabular} & Birational map to $\PP^4$ \\
\hline 
\hline 
 Quadric hypersurface in $\PP^5$ & no & all & projection from a point on it \\
\hline 
\rowcolor{yellow!25} Cubic hypersurface in $\PP^5$ & not known & \begin{tabular}{c} just many examples but \\ with a precise conjecture \end{tabular} & \begin{tabular}{c} some constructions in \\ \cite{Morin,Fano,Hassett,RS1,RS3} \end{tabular} \\
\hline 
\rowcolor{yellow!8} Quartic hypersurface in $\PP^5$ & yes \cite{Totaro2015} & no known examples &  \\
\hline 
 \begin{tabular}{c}Complete intersection of \\two  quadrics in $\PP^6$ \end{tabular} & no & all & projection from a line on it \\
\hline 
\rowcolor{yellow!8} \begin{tabular}{c} Complete intersection of \\ a quadric and a cubic in $\PP^6$ \end{tabular} & yes \cite{NOtte} & no known examples & \\
\hline 
 \rowcolor{yellow!8} \begin{tabular}{c} Complete intersection of \\ three quadrics in $\PP^7$ \end{tabular} & yes \cite{HPT2} & just some examples & \begin{tabular}{c} \emph{e.g.}, when the fourfold contains a plane \\ then one takes the projection from it  \end{tabular} \\
\hline 
 Del Pezzo fourfold $\mathbb{G}(1,4)\cap\PP^7$ & no & all & projection from the unique $\sigma_{2,2}$-plane \\
\hline 
\rowcolor{yellow!25} \begin{tabular}{c} Gushel-Mukai fourfold \end{tabular} & not known & \begin{tabular}{c} just the examples in \\ Table~\ref{tabella1}  \end{tabular} & \begin{tabular}{c} see Table~\ref{tabella1} \end{tabular} \\
\hline 
 \begin{tabular}{c} Linear section in $\PP^9$ \\ of the spinorial  $\mathbb{S}^{10}\subset\PP^{15}$ \end{tabular}  & no & all \cite{Roth1949} & projection from a tangent space \\
\hline 
 \begin{tabular}{c} Linear section in $\PP^{10}$ \\ of $\mathbb{G}(1,5)\subset\PP^{14}$ \end{tabular}  & no & all \cite{Roth1949} & \begin{tabular}{c} projection from the linear span of a \\ quintic del Pezzo surface contained in it  \end{tabular} \\
\hline 
 \begin{tabular}{c} Linear section in $\PP^{11}$  of the \\ Lagrangian Grass. 
$LG(3,6)\subset\PP^{13}$ \end{tabular} & no & all \cite{Roth1949} & \begin{tabular}{c} linear system of hyperplane sections \\ through a line and with one double point \end{tabular} \\
\hline 
 \begin{tabular}{c} Fourfold of degree $18$, \\ genus $10$, and coindex $3$ in $\PP^{12}$ \end{tabular} & no & all \cite{Roth1949} & \begin{tabular}{c} linear system of hyperplane sections \\ through a conic and with one double point \end{tabular} \\
\hline 
\end{tabular}
\end{adjustbox}
\caption{Rationality of smooth prime Fano fourfolds of coindex $\leq 3$, \cite{fujita-polarizedvarieties,mukai-biregularclassification}.}
\label{tabella2}
\end{table}

\clearpage


\providecommand{\bysame}{\leavevmode\hbox to3em{\hrulefill}\thinspace}
\providecommand{\MR}{\relax\ifhmode\unskip\space\fi MR }
\providecommand{\MRhref}[2]{%
  \href{http://www.ams.org/mathscinet-getitem?mr=#1}{#2}
}
\providecommand{\href}[2]{#2}

\end{document}